\documentclass[12pt,a4paper,reqno]{amsart} 
\usepackage{amssymb}
\usepackage{latexsym}
\usepackage{amsmath}
\usepackage{mathrsfs}
\usepackage{calc}                         

\newcommand{\scal}[2]{\langle #1,#2\rangle}

\newcommand{\rr}[1]{\mathbb R^{#1}}

\newcommand{\zz}[1]{\mathbb Z^{#1}}

\newcommand{\nm}[2]{\Vert #1\Vert _{#2}}

\newcommand{\ep}{\varepsilon}

\newcommand{\cdo}{\, \cdot \, }

\newcommand{\vrum}{\vspace{0.1cm}}

\newcommand{\maclS}{\mathcal S}

\newtheorem{thm}{Theorem}
\newtheorem{prop}[thm]{Proposition}
\newtheorem{cor}[thm]{Corollary}

\theoremstyle{definition}

\theoremstyle{remark}

\author{Carmen Fern\'andez}

\address{Departamento de An\'alisis Matem\'atico, Universitat de
Val\`encia, Valencia, Spain}

\email{fernand@uv.es}

\author{Antonio Galbis}

\address{Departamento de An\'alisis Matem\'atico, Universitat de
Val\`encia, Valencia, Spain}

\email{antonio.galbis@uv.es}

\author{Joachim Toft}

\address{Department of Mathematics,
Linn{\ae}us University, V{\"a}xj{\"o}, Sweden}

\email{joachim.toft@lnu.se}

\title{Convenient descriptions of weight functions in
time-frequency analysis}

\begin{document}

\begin{abstract}
Let $v$ be a submultiplicative weight. Then we prove that
$v$ satisfies GRS-condition, if and only if $v\cdot e^{-\ep |\cdo |}$ is bounded
for every positive $\ep$. We use this equivalence to establish identification
properties between weighted Lebesgue spaces, and between
certain modulation spaces and Gelfand-Shilov spaces.
\end{abstract}

\maketitle

\section{Introduction}\label{sec0}

\par

In \cite{GRS}, Gel'fand, Raikov and Shilov considered a family of matrix
classes in $\zz d\times \zz d$, where each matrix class, $\mathcal A_v$,
depends on a weight function $v$ on $\rr d$, and consists of all matrices
$A=\big ( a(j,k) \big )_{j,k\in \zz d}$ such that
$$
\sup _{j\in \zz d}a(j,j-k)v(k) \in \ell ^1.
$$
Here the weight $v$ is positive and submultiplicative, i.{\,}e. $v$ is
even and fulfills $v(x+y)\le v(x)v(y)$.

\par

Since \cite{baskakov,GRS}, several investigations have been performed which confirm
the importance of the class $\mathcal A_v$. For example if $v \ge c$ for
some positive constant $c$, then $\mathcal A_v \subseteq \mathcal
B(\ell ^2(\zz d))$, the set of all matrices which are bounded on $\ell ^2(\zz d)$.
Furthermore, in \cite{gr-book, gr_sjostrand}, strong links
between $\mathcal A_v$ and certain modulation spaces are described,
which appear in natural ways when considering non-stationary filters
in signal analysis (cf. \cite{strohmer} and the references therein).

\par

In \cite{baskakov,GRS}, the condition
\begin{equation}\label{GRScondition}
\lim _{\ell \to \infty}v(\ell x)^{1/\ell } =1,\quad x\in \rr d,
\end{equation}
is introduced for weights $v$ on $\rr d$, and link this condition to
fundamental questions on inverse closed properties for matrix algebras.
More precisely, for any submultiplicative weight $v$ which fulfills
\eqref{GRScondition}, it is proved in \cite{GRS} that the matrix
class $\mathcal A_v$ is so-called \emph{inverse closed} on $\ell ^2(\zz d)$.
That is, $\mathcal A_v$ is contained in $\mathcal B(\ell ^2(\zz d))$, and if
$A\in A_v$ is invertible on $\ell ^2(\zz d)$, then its inverse $A^{-1}$, does 
also
belong to $A_v$. Furthermore, if $v$ is submultiplicative and inverse closed
on $\ell ^2(\zz d)$, then it is necessary that $v$ satisfies 
\eqref{GRScondition}
(cf. \cite[Corollary 5.31]{gr-survey}).

\par

The condition \eqref{GRScondition} (the so-called
\emph{GRS-condition} or \emph{Gel'fand-Raikov-Shilov condition})
is therefore strongly linked to inverse closed properties for matrix algebras,
and a submultiplicative weight which satisfies the GRS-conditions is called
a \emph{GRS weight}.

\par

Later on, similar links between the GRS-conditions and inverse closed
properties for pseudo-differential operators with symbols in certain
modulation spaces were established (see e.{\,}g. \cite{gr_sjostrand}).

\par

In \cite{gr_fieldspaper}, it is proved that if $v$ is submultiplicative, then
there are positive constants $c$ and $C$ such that
$$
v(x)\le Ce^{c|x|},
$$
and it is easily seen that the right-hand side is a submultiplicative weight.
We note however that the weight
\begin{equation}\label{ExpWeight}
v_c(x) =v_{c,1}(x)\equiv e^{c|x|}
\end{equation}
does not satisfy the GRS-condition. Consequently, the matrix algebra
$\mathcal A_v$ is not inverse closed when $v=v_c$ and $c>0$.
On the other hand, it is easily seen that the set of submultiplicative
weights satisfying the GRS-condition is a monoid and contains
the weights
$$
x\mapsto (1+|x|)^r \quad \text{and}\quad x\mapsto
v_{c,s}(x)\equiv e^{c|x|^s}
$$
when $r\ge 0$,  $c\ge 0$ and $0\le s<1$.

\par

Since the (non-GRS) weight $v_{c,1}$ is rather similar to the GRS
weight $v_{c,s}$ when $s<1$ is close to $1$, it might be anticipated
that the the algebra $\mathcal A_{v_{c,1}}$ should possess a weaker
form of inverse closed property. Such property is established in
\cite{jaffard}, where it is proved that if $A\in \mathcal A_{v_{c_1,1}}$ for some
$c_1>0$ is invertible on $\ell ^2(\zz d)$, then its inverse
$A^{-1}$ belongs to $\in \mathcal A_{v_{c_2,1}}$, for some $c_2$.

\par

In this context we remark that in \cite{FGT}, a condition on sequences
of weights, which extends the notion of GRS-condition and includes exponential 
type weights is considered. Furthermore, certain type of
inverse closed property is established for matrix algebras which are
parameterized with such sequences of weights. In particular, the classical
equivalence between the inverse closed property of $\mathcal A_v$
and the fact that $v$ should be a GRS weight becomes a special
case of \cite[Theorem 2.1]{FGT}.

\medspace

In this paper we establish a more narrow link between the GRS-condition
and weights bounded by exponentials. More precisely, if $v$ is 
submultiplicative,
then we prove that $v$ fulfills the GRS-conditions, if and only if for every
$\ep >0$, there is a constant $C_\ep$ such that $v\le C_\ep v_{\ep ,1}$.

\par

We also combine this result with Theorem 3.9 in \cite{toft} to establish
relations with modulation spaces parameterized with weights satisfying
GRS-conditions, and Gelfand-Shilov spaces.

\par

\section{Identifications of the GRS-condition in terms of
exponentials}\label{sec1}

\par

We recall that a weight on $\rr d$ is a positive function on $\rr d$
which belongs to $L^\infty _{loc}(\rr d)$. The weight $v$ on $\rr d$
is submultiplicative, if $v$ is positive and even, and satisfies
$$
v(x+y)\leq v(x)v(y),\quad \mbox{for}\quad x,y \in \rr d.
$$
It is well-known that
$$
v(x)\lesssim e^{c|x|}
$$
for some $c>0$, when $v$ is submultiplicative. Furthermore, the
weight $v$ satisfies the \emph{GRS-condition} whenever
\begin{equation}\label{eq:grslimite0}
\lim _{\ell \to \infty}v(\ell x)^{1/\ell }=1,\quad \text{when}\quad x\in \rr d,
\end{equation}
holds.

\par

We have now the following observation concerning submultiplicative weights.

\par

\begin{prop}\label{prop:main}
Let $v$ be a submultiplicative weight on $\rr d$. Then the following conditions
are equivalent:
\begin{enumerate}
\item $v$ satisfies the GRS-condition;

\vrum

\item $v$ satisfies $v(x)\lesssim e^{\ep |x|}$, for every $\ep >0$.
\end{enumerate}
\end{prop}

\par

\begin{proof}
First assume that (2) holds. Let $\ep >0$ and $x\in \rr d$ be arbitrary and
fixed. As 
\begin{equation*}
v(nx) \le C e^{\ep n|x|},
\end{equation*}
we have $$
v(nx)^{1/n} \le  C^{1/n}e^{\ep |x|},
$$
and since $\ep $ is arbitrary we have
$$
\limsup _{n\to \infty} v(nx)^{1/n} \le 1.
$$

\par

On the other hand, since $v(x)\ge c$ for some $c>0$, it follows that
$$
\liminf _{n\to \infty} v(nx)^{1/n} \ge 1.
$$
A combination of these estimates gives that $\lim _{n\to \infty} v(nx)^{1/n}$ exists
and is equal to $1$, i.{\,}e. (1) holds.

\par

Next we assume that (2) does not hold. Let $\{e_1,\dots,e_d\}$
be the canonical  basis  in $\mathbb{R}^d.$ Since $v$ is
submultiplicative,
$$
v(x_1,\dots,x_d)\leq  \prod_{j=1}^{d}v(x_je_j).
$$
Hence, if (2) is not fulfilled, then there are $j\in \{1, 
\dots, d\}$ and $\ep >0$ such that $v(te_j)e^{-\ep |t|}$ is
unbounded. Without loss of generality we may assume that $j=1$. The weight on $\mathbb{R}$ defined by 
$$
u(t):=v(te_1)
$$
is submultiplicative. Since $u(t)e^{-\ep t}$ is 
unbounded in $[0,\infty[$, we may find a sequence $\left(t_n\right)$
such that  $0<t_n+1<t_{n+1}\to \infty$ and
$$
u(t_n) > e^{\ep t_n}.
$$
Let $k_n$ be the 
integer part of $t_n.$ As $t_n<k_n+1$ and $u$ is submultiplicative,
$$
u(t_n)\le  C u(k_n+1)\leq C u(1)u(k_n),
$$
where $C = \max_{|\xi|\leq 1}u(\xi)$. Hence 
$$
u(k_n)\geq \frac{u(t_n)}{C u(1)}\geq \frac{e^{\ep t_n}}{C u(1)}
$$
and 
$$
u(k_n)^{1/k_n}\geq \left(\frac{u(t_n)}{C u(1)}\right)^{1/t_n}\geq \frac{e^{\ep 
}}{(C u(1))^{1/t_n}}.
$$
Consequently
$$
\liminf_{n\to \infty} u(k_n)^{1/k_n} \geq 
e^{\ep }.
$$
This implies that $u,$ and therefore $v,$ does not satisfy the 
GRS-condition. The proof is complete.
\end{proof}

\par\medskip\
The previous result means that ${\mathcal A}_v$ is inverse closed in $\mathcal B(\ell ^2(\zz d))$ if, and only if, it contains the Jaffard algebra of all matrices with exponential off-diagonal decay.

\begin{cor}\label{CorProp}
Let $\Omega$ be the set of all weights on $\rr {2d}$ which satisty
the GRS-condition. Then
\begin{equation}\label{ellAidentities}
\bigcup_{\ep >0}\ell^1_{(v_\ep )} =\bigcap_{v\in \Omega}
\ell^1_{(v)}
\quad \text{and} \quad
\bigcup_{\ep >0}\mathcal A_{v_\ep}
=\bigcap_{v\in \Omega} \mathcal A_v.
\end{equation}
\end{cor}

\par

\begin{proof}
This is a consequence of Proposition \ref{prop:main} and 
\cite[Theorem 2.1]{FGT}.
\end{proof}

\par

Next we analyze a condition on sequences of weights
considered in  \cite{FGT}.

\par

\begin{prop}
Let ${\mathcal W} = \{ w_n\} _{n\in \mathbb N}$ be a decreasing sequence of 
submultiplicative weights. Then the following conditions are
equivalent:
\begin{itemize}
\item[(1)] $\begin{displaystyle}\inf_n \lim_{\ell \to \infty}w_n(\ell 
x)^{{1}/{\ell}} = 1\end{displaystyle}$ for every $x\in {\mathbb R}^d$.

\vrum

\item[(2)] for every $\ep > 0$, there exists an $n\in {\mathbb N}$ such
that $w_n(x)e^{-\ep |x|}$ is bounded.
\end{itemize}
\end{prop}

\par

\begin{proof}
$(2)\Rightarrow (1)$ is clear. Assume that (1) holds, and let
$V$ be the family of all submultiplicative 
and weights $v\, :\,  \zz d\to \mathbb R_+$ such that 
$$
\sup_{{k}}\frac{v(k)}{w_n(k)} < \infty\quad \forall \ n\in
{\mathbb N}.
$$
By \cite[Theorem 2.3]{bms} we have
$$
\bigcap_{v\in V}\ell^1_{(v)} = \bigcup_n \ell^1_{(w_n)}
$$
algebraic and topologically (see also the proof of \cite[Theorem 2.1]{FGT}).
Since every weight $v\in V$ satisfies the GRS-condition, we have $v = 
O(v_{\ep})$ for every $\ep > 0.$ 
Consequently,
$$
\ell^1 _{(v_{\ep})} \subset \bigcup_n \ell^1_{(w_n)}
$$
with continuous inclusion. The result now follows from
Gr\"othendieck's  factorization theorem (see \cite[1.2.20]{bp}).
\end{proof}

\par

We finish the paper by applying our previous results to establish certain
links between modulation spaces and Gelfand-Shilov spaces.
We recall the definition of these spaces.

\par

Let $0<h,s\in
\mathbb R$ be fixed, and let $\mathcal S_{s,h}(\rr d)$ be the set of
all $f\in C^\infty (\rr d)$ such that
\begin{equation*}
\nm f{\mathcal S_{s,h}}\equiv \sup \frac {|x^\beta \partial ^\alpha
f(x)|}{h^{|\alpha | + |\beta |}\alpha !^s\, \beta !^s}
\end{equation*}
is finite. Here the supremum should be taken over all $\alpha ,\beta \in
\mathbb N_0^d$ and $x\in \rr d$. We note that $\mathcal S_{s,h}(\rr d)
\subseteq \mathscr S(\rr d)$, increases with $h$ and $s$, where
$\mathscr S(\rr d)$ is the set of all Schwartz functions on $\rr d$.
The \emph{Gelfand-Shilov space} $\mathcal S_{s}(\rr d)$ is the inductive limit
of $\mathcal S_{s,h}(\rr d)$ and the
\emph{Gelfand-Shilov distribution space} $\mathcal S_{s}'(\rr d)$
is the projective limit of $\mathcal S_{s,h}'(\rr d)$.
In particular,
\begin{equation}\label{GSspacecond1}
\mathcal S_{s}(\rr d) = \bigcup _{h>0}\mathcal S_{s,h}(\rr d),
\quad \text{and}\quad
\mathcal S_s'(\rr d) = \bigcap _{h>0}\mathcal S_{s,h}'(\rr d).
\end{equation}
Note that $\mathcal S_s'(\rr d)$
is the dual of $\mathcal S_s(\rr d)$ (\cite{GS, Ko, Pil}).

\par

Evidently, $\mathscr S(\rr d)\subset \mathcal S_s(\rr d)$ increases with
$s$, and $\mathscr S'(\rr d)\subseteq \mathcal S_s'(\rr d)$
decreases with $s$. If $s<1/2$, then $\mathcal S_s(\rr d)$ is trivial.

\par

The Fourier transform on $\mathscr S'(\rr d)$ is the linear and
continuous map which takes the form
$$
(\mathscr Ff)(\xi )= \widehat f(\xi ) \equiv (2\pi )^{-d/2}\int _{\rr
{d}} f(x)e^{-i\scal  x\xi }\, dx
$$
when $f\in L^1(\rr d)$. For every $s\ge 1/2$, the Fourier transform
is continuous and bijective on $\mathcal S_s(\rr d)$, and extends
uniquely to a continuous and bijective map on $\mathcal S_s'(\rr d)$.

\par

Next we recall the definition of modulation spaces. Let $1/2\le s_0<s$,
and let $\phi \in \maclS _{s_0}(\rr d)\setminus 0$ be fixed.
Then the \emph{short-time Fourier transform} $V_\phi f$ of
$f\in \maclS _s'(\rr d)$ is the element in $\maclS _s'(\rr {2d})
\cap C^\infty (\rr {2d})$, defined by the formula
$$
(V_\phi f)(x,\xi ) := \mathscr F(f\cdot \phi (\cdo -x))(\xi).
$$
If in addition $f$ is in the Schwartz class, then
$V_\phi f$ is given by
$$
(V_\phi f)(x,\xi ) = (2\pi )^{-d/2}\int _{\rr d} f(y)\overline{\phi (y-x)}
e^{-i\scal y\xi}\, dy.
$$

\par

Let $v$ be a submultiplicative weight on $\rr d$. Then the weight
$\omega$ on $\rr d$ is called $v$-moderate if
$$
\omega (x+y)\le C\omega (x)v(y),
$$
for some constant $C$ which is independent of $x\in \rr d$ and
$y\in \rr d$.

\par

Let $\phi \in \maclS _{1/2}(\rr d)$, $1/2<s<1$, $p,q\in (0,\infty]$, $v$ be
submultiplicative on $\rr {2d}$, and let $\omega$ be a $v$-moderate
weight on $\rr {2d}$. Then the \emph{modulation space}
$M^{p,q}_{(\omega )}(\rr d)$ is the quasi-Banach space which consists of
all $f\in \maclS _{1/2} '(\rr d)$ such that
\begin{equation}\label{modnorm}
\nm f{M^{p,q}_{(\omega )}} := \left ( \int _{\rr d} \left ( \int _{\rr d}
|V_\phi f (x,\xi )\omega (x,\xi ) |^p\, dx \right )^{q/p}\, d\xi \right )^{1/q}
\end{equation}
is finite (with obvious modifications when $p=\infty$ or $q=\infty$).
The definition of $M^{p,q}_{(\omega )}(\rr d)$ is
independent of the choice of $\phi \in \maclS _{1/2} (\rr d)\setminus 0$,
and different $\phi$ gives rise to equivalent norms. (See e.{\,}g.
\cite{gr-book}).

\medspace

According to \cite[Prop. 11.3.1]{gr-book}, 
$$
\bigcap_{s\geq 0}M^{\infty}_{(w_s)}({\mathbb R}^d) =
{\mathscr S}({\mathbb R}^d),
$$
where $w_s(x) = (1+|x|)^s$. However, when considering all the weights
satisfying the GRS-condition, the intersection is no longer a Fr\'echet
space but an inductive limit of Banach spaces.
%
%

\par

\begin{prop}
Let $p,q\in (0,\infty ]$, and let $\Omega$ be the set of all submultiplicative
weights on $\rr {2d}$ which satisty the GRS-condition.. Then
\begin{equation}\label{capcupidentities}
\bigcap_{v\in \Omega} M^{p,q}_{(v)}(\rr d) = \mathcal S_1(\rr d)
\quad \text{and}\quad
\bigcup_{v\in \Omega} M^{p,q}_{(1/v)}(\rr d) = \mathcal S_1'(\rr d).
\end{equation}
\end{prop}

\par

For the proof we note that if $\omega$ is a moderate weight on
$\rr {2d}$, $p\in (0,\infty ]$ and $f\in \maclS '_1(\rr d)$, then $f\in
M^{p,p}_{(\omega )}(\rr d)$, if and only if
\begin{equation}\label{CubeCoeff}
\left \{   \left ( \iint_{n+Q} |V_\phi f(x,\xi )|^p\, dxd\xi 
\right )^{1/p} \right \} _{n\in \zz {2d}} \in
\ell ^{p}_{(\omega )}(\zz {2d}).
\end{equation}
Here $Q$ is the cube $[0,1]^{2d}$.

\par

\begin{proof}
By well-known embedding properties for modulation spaces, it
suffices to prove the result in the case $p=q=1$ or $p=q=\infty$
(cf. \cite{gr-book}). Let $\phi \in
{\mathcal S}_{1/2}({\mathbb  R}^d)\setminus \left\{0\right\}$ be fixed
and let $Q:=[0,1]^{2d}$ be the unit  cube. For any
submultiplicative weight $v$ on $\rr {2d}$ and $f\in
{\mathcal S}^\prime _1({\mathbb R}^d)$ it turns out that 
$f\in M^1_{(v)}({\mathbb R}^d)$, if and only if \eqref{CubeCoeff}
holds for $p=1$. Consequently, the first identity in \eqref{ellAidentities} 
gives 
$$
\bigcap_{v\in \Omega} M^{1}_{(v)}(\rr d) = \bigcup_{\epsilon > 0} M^{1}
_{(v_{\epsilon})}(\rr d) = \mathcal S_1(\rr d),
$$
where the last identity follows from Theorem 3.9 in 
\cite{toft}. The second statement
of the proposition follows from the fact that the dual of $M^1_{(v)}$
is given by $M^\infty _{(1/v)}$ (see e.{\,}g.  \cite{gr-book}).
\end{proof}

\vskip 0.5cm
\textbf{Acknowledgement.} The research of the first two authors was partially
supported by MEC and FEDER Project MTM2010-15200 and GVA
Prometeo II/2013/013.

\end{document}